\documentclass[a4paper,12pt,twoside]{article}
\usepackage{amsfonts,theorem}

\newcommand{\C}{\mathbb{C}}

\renewcommand{\P}{\mathbb{P}}

\newcommand{\R}{\mathbb{R}}

\renewcommand{\L}{\mathcal{L}}
\renewcommand{\O}{\mathcal{O}}

\renewcommand{\a}{\mathfrak{a}}

\renewcommand{\l}{\mathfrak{l}}
\newcommand{\suchthat}{\mathop{\,\vert\,}}
\makeatletter
\def\operatorname#1{\mathop{\operator@font #1}\nolimits}%
\makeatother

\newcommand{\End}{\operatorname{End}}

\newcommand{\Id}{\operatorname{Id}}

\newcommand{\Tr}{\operatorname{Trace}}
\renewcommand{\div}{\operatorname{div}}

\newtheorem{thm}{Theorem}

\newtheorem{cor}{Corollary}
\newtheorem{lem}{Lemma}
{\theorembodyfont{\normalfont\rmfamily}

\newtheorem{rem}{Remark}
}
\makeatletter
\newenvironment{pf}[1][]%
 {\def\proof@temp{#1}\par\noindent
  \textsc{Proof}\ifx\proof@temp\@empty\else\ (#1)\fi\hspace{1em}}
 {\rule{12pt}{0pt}\hfill\rule{6pt}{6pt}\par\vspace{.4\baselineskip}}
\makeatother

\makeatletter
\renewcommand{\thesection}{\@arabic\c@section}
\def\section{\@startsection {section}{1}{0pt}{12pt plus 2pt
minus 2pt}{-4pt}{\bf\setcounter{equation}{0}}}
\makeatother

\makeatletter
\def\footnotenomark#1{
    \protected@xdef\@thanks{\@thanks
        \protect\footnotetext{#1}}%
}
\makeatother
\newcommand{\note}[1]{${}^{\hbox{\small(#1)}}$}

\title{Homogeneous symplectic manifolds \\ 
with Ricci-type curvature}

\author{M.~Cahen\note{i},\ S.~Gutt\note{i,ii}\footnotenomark{Research of
the first three authors supported by an ARC of the communaut\'e fran\c caise
de Belgique.}\footnotenomark{The third author is Aspirant au Fonds National
 belge de la Recherche Scientifique.}\footnotenomark{\textbf{Mathematics Subject 
Classification (1991):}\ \ 53C05, 58C35, 53C57.},\\[5pt] 
J.~Horowitz\note{i} and J.~Rawnsley\note{ii,iii}
\footnotenomark{Email: mcahen@ulb.ac.be, sgutt@ulb.ac.be, 
horowitz@ulb.ac.be and j.rawnsley@warwick.ac.uk}\\[30pt]
\small\makebox[.5in][r]{(i)}
\parbox[t]{2.5in}{\noindent Universit\'e Libre de Bruxelles\\
Campus Plaine, CP 218\\
bvd du triomphe\\
1050 Brussels\\
Belgium}\\~\\[2pt]
\small\makebox[.5in][r]{(ii)}
\parbox[t]{2.5in}{\noindent Universit\'e de Metz\\
Ile du Saulcy\\
57045 Metz Cedex 01\\
France}\\~\\[2pt]
\small\makebox[.5in][r]{(iii)}
\parbox[t]{2.5in}{\noindent Mathematics Institute\\
University of Warwick\\
Coventry CV4 7AL\\
United Kingdom}\\[10pt]~}

\date{\small May 2000\\[-30pt]~}

\renewcommand{\baselinestretch}{1.25}

\begin{document}

\setcounter{page}{0}
\renewcommand{\baselinestretch}{1} 
\thispagestyle{empty}
\maketitle

\begin{abstract}
We consider invariant symplectic connections $\nabla$ on homogeneous
symplectic manifolds $(M,\omega)$ with curvature of Ricci type. Such
connections are solutions of a variational problem studied by Bourgeois
and Cahen, and provide an integrable almost complex structure on the
bundle of almost complex structures compatible with the symplectic
structure. If $M$ is compact with finite fundamental group then
$(M,\omega)$ is symplectomorphic to $\P_n(\C)$ with a multiple of its
K\"ahler form and $\nabla$ is affinely equivalent to the Levi-Civita
connection.
\end{abstract}

\thispagestyle{empty}

\newpage

\renewcommand{\baselinestretch}{1.25}

The space of curvature tensors of symplectic connections on a symplectic
manifold $(M,\omega)$ of dimension $2n\ge4$ splits under the action of
the symplectic group $Sp(2n,\R)$ as a direct sum of two subspaces on
which $Sp(2n,\R)$ acts irreducibly
\cite{bib:BourgCah,bib:Vaisman1,bib:MaudDV}. For a given
curvature tensor $R$ we shall denote by $E$ and $W$ its projections onto
these two subspaces. The $E$-component is determined by the Ricci
tensor of the connection. When the $W$-component vanishes identically we
say that the curvature is of Ricci type.

The motivation for looking at such connections is two-fold. They provide
critical points of a functional which has been introduced in
\cite{bib:BourgCah} to select preferred symplectic connections, and
$W=0$ is the integrability condition for an almost complex structure
which a symplectic connection determines on the total space of the
bundle $J(M,\omega)$ of almost complex structures compatible with the
symplectic structure \cite{bib:OBR,bib:Vaisman2,bib:Vaisman3}.

The simplest framework in which one can study the $W=0$ condition is the
compact homogeneous one. Our main result is

\begin{thm}\label{thm:1}
Let $(M,\omega)$ be a compact homogeneous symplectic manifold with
finite fundamental group. If $(M,\omega)$ admits a homogeneous
symplectic connection $\nabla$ with Ricci-type curvature then
$(M,\omega)$ is symplectomorphic to $(\P_n(\C),\omega_0)$, where
$\omega_0$ is a multiple of the K\"ahler form of the Fubini--Study
metric, and $\nabla$ is affinely equivalent to the Levi-Civita
connection.
\end{thm}

When we do not impose any restriction on the
fundamental group, we were only able to prove

\begin{thm}\label{thm:2}
Let $(M,\omega)$ be a compact homogeneous symplectic manifold of
dimension $4$. If $(M,\omega)$ admits a homogeneous symplectic
connection $\nabla$ with Ricci-type curvature then $\nabla$ is locally
symmetric.
\end{thm}

In $\S$\ref{section:1} we prove some general identities which hold for
any symplectic connection with Ricci-type curvature. In
$\S$\ref{section:2} we deduce some easy consequences of these
identities in the homogeneous (respectively compact homogeneous)
framework. In $\S$\ref{section:3} we prove Theorem \ref{thm:1}  in
the simply connected case and show how to extend this to a finite
fundamental group. Finally $\S$\ref{section:4} is devoted to the proof
of Theorem \ref{thm:2}.

\section{} \label{section:1} Let $(M,\omega)$ be a symplectic manifold
and $\nabla$ be a symplectic connection (a torsion-free connection on
$TM$ with $\nabla\omega=0$). The curvature endomorphism $R$ of $\nabla$
is defined by
\[
R(X,Y)Z = \left(\nabla_X\nabla_Y - \nabla_Y\nabla_X -
\nabla_{[X,Y]}\right)Z
\]
for vector fields $X,Y,Z$ on $M$. The symplectic curvature tensor
\[
R(X,Y;Z,T) = \omega(R(X,Y)Z,T)
\]
is antisymmetric in its first two arguments, symmetric in its last two
and satisfies the first Bianchi identity
\[
\oint_{X,Y,Z} R(X,Y;Z,T) = 0
\]
where $\oint$ denotes the sum over the cyclic permutations of the listed
set of elements. The Ricci tensor $r$ is the symmetric $2$-tensor
\[
r(X,Y) = \Tr[ Z \mapsto R(X,Z)Y].
\]
$R$ also obeys the second Bianchi identity
\[
\oint_{X,Y,Z} \left(\nabla_XR\right)(Y,Z) = 0.
\]
The Ricci part $E$ of the curvature tensor is given by
\begin{eqnarray}\label{one:1.0}
E(X,Y;Z,T) &=& \frac{-1}{2(n+1)} \biggl[2\omega(X,Y)r(Z,T) +
\omega(X,Z)r(Y,T) + \omega(X,T)r(Y,Z)\nonumber\\
&&\qquad\mbox{}-\omega(Y,Z)r(X,T) - \omega(Y,T)r(X,Z)
\biggr].
\end{eqnarray}
The curvature is of Ricci type when $R=E$.

\begin{lem}
Let $(M,\omega)$ be a symplectic manifold of dimension $2n\ge4$. If the
curvature of a symplectic connection $\nabla$ on $M$ is of Ricci type
then there is a $1$-form $u$ such that
\begin{equation}\label{one:1.1}
\left(\nabla_Xr\right)(Y,Z) = \frac{1}{2n+1}\left(\omega(X,Y)u(Z)
+ \omega(X,Z)u(Y)\right).
\end{equation}
Conversely, if there is such a $1$-form $u$ then the $W$ part of the
curvature satisfies
\begin{equation}
\oint_{X,Y,Z} \left(\nabla_X W\right)(Y,Z;T,U) = 0.
\end{equation}
\end{lem}

\begin{pf}
When the curvature is of Ricci type, the second Bianchi identity for $R$
becomes an identity for $E$. Since $\omega$ is parallel,
covariantly differentiating equation (\ref{one:1.0}) and summing cyclically,
we get
\begin{eqnarray}\label{one:*}
0 &=& \oint_{X,Y,Z} 2\omega(Y,Z)(\nabla_Xr)(T,U)
+\omega(Y,T)(\nabla_Xr)(Z,U)
+\omega(Y,U)(\nabla_Xr)(Z,T)\nonumber\\
&&\qquad\qquad\mbox{}-\omega(Z,T)(\nabla_Xr)(Y,U)
-\omega(Z,U)(\nabla_Xr)(Y,T).
\end{eqnarray}
Choose local frames $\{V_a\}_{a=1}^{2n}$, $\{W_a\}_{a=1}^{2n}$ on $M$
such that $\omega(V_a,W_b)=\delta_{ab}$. Substitute $Y=V_a$ and $Z=W_a$
in equation (\ref{one:*}) and sum over $a$ to obtain
\begin{eqnarray}\label{one:**}
0 &=&2n(\nabla_Xr)(T,U)
- (\nabla_Tr)(X,U)
- (\nabla_Ur)(X,T)\nonumber\\
&&\qquad\mbox{}+\omega(X,T)\sum_a(\nabla_{W_a}r)(V_a,U)
+\omega(X,U)\sum_a(\nabla_{W_a}r)(V_a,T).
\end{eqnarray}
If we cyclically permute $X,T,U$ in equation (\ref{one:**}) and sum we
get
\begin{equation}\label{one:***}
(2n-2) \oint_{X,T,U} (\nabla_Xr)(T,U) = 0
\end{equation}
and since $n\ge2$ we have
\begin{equation}\label{one:****}
\oint_{X,T,U} (\nabla_Xr)(T,U) = 0
\end{equation}
Using equation (\ref{one:****}) in equation (\ref{one:**}) gives
\[
(2n+1)(\nabla_Xr)(T,U) + \omega(X,T)\sum_a(\nabla_{W_a}r)(V_a,U)
 + \omega(X,U)\sum_a(\nabla_{W_a}r)(V_a,T) = 0
\]
which is of the desired form if
\[
u(X) = -\sum_a(\nabla_{W_a}r)(V_a,X).
\]

Conversely, if one substitutes (\ref{one:1.1}) into the covariant
derivative of (\ref{one:1.0}) and cyclically sums then one obtains
\[
\oint_{X,Y,Z} (\nabla_XE)(Y,Z,T,U)=0.
\]
Combining this with the second Bianchi identity, gives the second part
of the Lemma.
\end{pf}

\begin{cor}
A symplectic manifold with a symplectic connection whose curvature is of
Ricci type is locally symmetric if and only if the $1$-form $u$, defined
in the Lemma, vanishes.
\end{cor}

\begin{rem}
It will be useful to have an equivalent form of formula (\ref{one:1.1}).
Denote by $A$ the linear endomorphism such that
\begin{equation}\label{one:1.3}
r(X,Y) = \omega(X, AY).
\end{equation}
The symmetry of $r$ is equivalent to saying that
$A$ is in the Lie algebra of the symplectic group of $\omega$. Denote by
$\overline{u}$ the vector field such that
\begin{equation}\label{one:1.4}
u = i(\overline{u})\omega
\end{equation}
then (\ref{one:1.1}) is equivalent to
\begin{equation}\label{one:1.5}
\nabla_XA = \frac{-1}{2n+1}(X\otimes u + \overline{u}\otimes i(X)\omega).
\end{equation}
\end{rem}

\begin{lem}
Let $(M,\omega)$ be a symplectic manifold with a symplectic connection
$\nabla$ with Ricci-type curvature. Then,
keeping the above notation, the following identities hold:
\begin{enumerate}
\item There is a function $b$ such that
\begin{equation}\label{one:1.6}
\nabla u = -\frac{1+2n}{2(1+n)}\stackrel{(2)}{r} + b\omega
\end{equation}
where $\stackrel{(2)}{r}$ is the $2$-form
\begin{equation}\label{one:1.7}
\stackrel{(2)}{r}(X,Y) = \omega(X, A^2Y).
\end{equation}
\item The differential of the function $b$ is given by
\begin{equation}\label{one:1.8}
db = \frac{1}{1+n} i(\overline{u})r.
\end{equation}
\item The covariant differential of $db$ is given by
\begin{equation}\label{one:1.9}
\nabla db = \frac{1}{1+n} \left[ -\frac{1}{1+2n} u\otimes u -
\frac{1+2n}{2(1+n)}\stackrel{(3)}{r} + br\right].
\end{equation}
where
\begin{equation}\label{one:1.10}
\stackrel{(3)}{r}(X,Y) = \omega(X,A^3Y).
\end{equation}
\end{enumerate}
\end{lem}

\begin{pf}
We can compute the action of the curvature on endomorphisms in two
different ways. On the one hand it is
\begin{eqnarray*}
R(X,Y)\cdot A &=& [R(X,Y),A]\\
 &=& R(X,Y)A - AR(X,Y)\\
&=&  -\frac{1}{2(n+1)}[ X\otimes \omega(A^2Y,\,.\,) - Y \otimes
\omega(A^2X,\,.\,)\\
&&\qquad\mbox{} +A^2Y \otimes \omega(X,\,.\,) - A^2X\otimes\omega(Y,\,.\,)].
\end{eqnarray*}
On the other hand the curvature is of Ricci type so that (\ref{one:1.5})
gives
\[
R(X,Y)\cdot A = \frac{1}{2n+1}[X\otimes \nabla_Yu - Y \otimes \nabla_X
u + \nabla_Y \overline{u} \otimes \omega(X,\,.\,) -
\nabla_X\overline{u}\otimes\omega(Y,\,.\,)].
\]
If we define an endomorphism $B$ of $TM$ by
\[
BY = \frac{2n+1}{2(n+1)} A^2Y + \nabla_Y \overline{u}
\]
then equality of the two right hand sides yields
\[
X \otimes \omega(BY, \,.\,) - Y \otimes \omega(BX,\,.\,) + BY \otimes
\omega(X,\,.\,) - BX \otimes \omega(Y,\,.\,) = 0
\]
whose only solution is
\[
B= b \Id.
\]
This gives
\[
\nabla_Y u = -\frac{2n+1}{2(n+1)} \omega(A^2Y, \,.\,) + b \omega(Y, \,.\,)
\]
which is equation (\ref{one:1.6}).

Antisymmetrising  (\ref{one:1.6}) we get
\[
du = -\frac{2n+1}{n+1} \stackrel{(2)}{r} + 2 b \omega.
\]
Taking the exterior derivative gives
\[
0 = -\frac{2n+1}{n+1} d\stackrel{(2)}{r} + 2 db \wedge\omega.
\]
But
\begin{eqnarray*}
d\stackrel{(2)}{r}(X,Y,Z) &=& \oint_{X,Y,Z} \omega(\nabla_XA^2Y,Z)\\
&=& -\frac{1}{2n+1}\oint_{X,Y,Z} \omega(u(AY)X +
\omega(X,AY)\overline{u},Z)\\
&&\qquad\mbox{}+\omega(u(Y)AX + \omega(X,Y)A\overline{u},Z)\\
&=& \frac{2}{2n+1} \oint_{X,Y,Z} \omega(X,Y)r(\overline{u},Z).
\end{eqnarray*}
Substituting,
\[
\left[ -\frac{1}{n+1}r(\overline{u},\,.\,) + db\right]\wedge\omega = 0.
\]
and in dimension $4$ or higher this implies
\[
db=  \frac{1}{n+1}r(\overline{u},\,.\,)
\]
which is (\ref{one:1.8}). Covariantly differentiating
\begin{eqnarray*}
(\nabla_Xdb)(Y) &=& \frac{1}{n+1}\left[(\nabla_Xr)(\overline{u},Y)
+ r(\nabla_X \overline{u},Y)\right]\\
&=&\frac{1}{n+1}\left[\frac{1}{1+2n}\omega(X,\overline{u})u(Y) +
r\left(-\frac{2n+1}{2(n+1)}A^2X+bX,Y\right)\right]
\end{eqnarray*}
which is (\ref{one:1.9}).
\end{pf}

\section{} \label{section:2}
Assume $(M,\omega)$ is a $G$-homogeneous symplectic manifold and
$\nabla$ is a $G$-invariant symplectic connection with Ricci-type
curvature. If $\nabla$ is not locally symmetric the $G$-invariant
$1$-form $u$ is everywhere different from zero and the function $b$ is also
$G$-invariant and hence constant. Putting these two facts into
(\ref{one:1.8}) we see that $r$ as a bilinear form is necessarily
degenerate
\begin{equation}\label{two:2.1}
r(\overline{u},\,.\,) = 0.
\end{equation}
Also (\ref{one:1.9}) implies
\begin{equation}\label{two:2.2}
\frac{1}{2n+1} u\otimes u + \frac{1+2n}{2(1+n)} \stackrel{(3)}{r} - br
=0
\end{equation}
or equivalently
\begin{equation}\label{two:2.3}
\frac{1}{2n+1} \overline{u}\otimes u - \frac{1+2n}{2(1+n)} A^3 + bA =0.
\end{equation}
Applying $A$ to (\ref{two:2.3}) and using (\ref{two:2.1})
\begin{equation}\label{two:2.4}
-\frac{1+2n}{2(1+n)} A^4 + bA^2 =0.
\end{equation}
It follows that the only possible non-zero eigenvalues of $A$ are
$\pm\sqrt{\frac{2(1+n)}{1+2n}b}$ and so are real or imaginary.

\begin{lem}\label{two:lem1}
If $(M,\omega)$ is a compact homogeneous symplectic manifold admitting a
homogeneous symplectic connection $\nabla$ with Ricci-type curvature
which is not locally symmetric then $b=0$.
\end{lem}

\begin{pf}
Recall that for any vector field $X$, Cartan's identity  gives
\[
\div X \omega^n \stackrel{\hbox{\small def}}{=\!=} \L_X \omega^n
= n\, d\left((i(X)\omega)\wedge\omega^{n-1}\right)
\]
and
\[
\L_X \omega^n = \left(\L_X  - \nabla_X\right)\omega^n =n\left(
\omega(\nabla_. X,\,.\,) + \omega(\,.\,,\nabla_. X)
\right)\wedge\omega^{n-1}
\]
so that
\[
\div X = \Tr[ Z \mapsto \nabla_ZX].
\]
In particular, by (\ref{one:1.6})
\[
\div \overline{u} = -\frac{2n+1}{2(n+1)} \Tr A^2 + 2nb.
\]
$G$-invariance implies that $\div \overline{u}$ is constant. But $M$
compact with no boundary implies $\int_M \div \overline{u} \omega^n = 0$
since the
argument is exact; hence the constant is zero. Thus
\[
b = \frac{2n+1}{4n(n+1)} \Tr A^2.
\]
On the other hand, (\ref{two:2.4}) implies that $A^2$ is a multiple of a
projection and with $A$ symplectic this has even rank $2p$ say; using
\ref{two:2.1} we get $2p<2n$. Thus
\[
\Tr A^2 = \frac{4pb(1+n)}{1+2n}
\]
so
\[
b = \frac{2n+1}{4n(n+1)}\,.\, \frac{4pb(1+n)}{1+2n} = \frac{p}{n} b
\]
and hence $b=0$.
\end{pf}

It follows that $A^4=0$ so $A$ is nilpotent; moreover (\ref{two:2.3})
tells us that $A^3$ has rank 1.

\begin{lem}\label{two:lem2}
Let $(M,\omega)$ be a $4$-dimensional homogeneous symplectic manifold
admitting a
homogeneous symplectic connection $\nabla$ with Ricci-type curvature
which is not locally symmetric. Let $A$ be the endomorphism associated
to the Ricci tensor. Then
\begin{enumerate}
\item either $A$ is nilpotent, $b \ne 0$, $A^2=0$,  and $A$ has rank 1 at
any point;
\item or  $A$ is nilpotent, $b = 0$, and $A^3$ has rank 1 at any point;
\item or $A$ has a non zero eigenvalue so $b\ne 0$. Then $A$
admits a pair of non zero eigenvalues of opposite sign (real or
imaginary) with multiplicity $1$ and $0$ is an eigenvalue of
multiplicity $2$ at any point. Furthermore, $A$ has necessarily a nilpotent
part.
\end{enumerate}
\end{lem}

\begin{pf}
The dimension -- at any point $x\in M$ -- of the  generalised $0$ eigenspace
of $A$ is even and non-zero, so
is $2$ or $4$. If it is $4$ then $A$ is nilpotent and $A^4=0$ in
dimension $4$. Thus, by (\ref{two:2.3}), $bA^2=0$. If $b\ne0$ then $A^2=0$ so,
 by (\ref{two:2.4}), $A$ has rank 1 at any point.
Otherwise $b=0$ and $A^3$ has rank 1 at any point.

When the generalised 0 eigenspace $V_0$
is $2$-dimensional at any point, then $\pm\sqrt{\frac{2(1+n)}{1+2n}b}$ are
eigenvalues
with multiplicity $1$.
Choose a globally defined vector field $v\in V_0$ so that
$\omega(v,\overline{u})=1$. Set $Av=p\overline{u}$. Then
\[
\nabla_X(Av)=\nabla_X(p\overline{u})=(Xp)\overline{u}
+ p\left( -\frac{1+2n}{2(1+n)}A^2X+bX\right)
\]
but it is also equal to
\[
\nabla_X(Av)=(\nabla_XA)v+A(\nabla_Xv)=-\frac{1}{1+2n}\left(X
u(v)+\overline{u}\omega(X,v)\right)
+A(\nabla_Xv).
\]
Observe that $\omega(A^2X,\overline{u})=\omega(A(\nabla_Xv),\overline{u})=0$,
so that
\[
p\omega(bX,\overline{u})=\omega(-\frac{1}{1+2n}u(v)X ,\overline{u})=
\frac{1}{1+2n}\omega(X ,\overline{u}).
\]
Hence $pb=\frac{1}{1+2n}$ which implies that $p\ne 0$. Thus $A$ has a nilpotent
part.
\end{pf}

\section{} \label{section:3}
We first prove Theorem \ref{thm:1} in the simply-connected case. It is
standard that a compact simply-connected homogeneous symplectic manifold
$(M,\omega)$ is symplectomorphic to a coadjoint orbit of a
simply-connected compact semisimple Lie group $G$. Such a Lie group $G$
is a product of simple groups and the orbit is a product of orbits. We
may throw away any factors where the orbit is zero dimensional as the
remaining group will still act transitively. A $G$-invariant symplectic
connection $\nabla$ on such an orbit is compatible with the product
structure. If the curvature of $\nabla$ is of Ricci type, then it was
shown in \cite{bib:CGR2}  that the curvature is zero when
$(M,\omega,\nabla)$ is a product of more than one factor. But a
non-trivial compact coadjoint orbit of a simple Lie group does not admit
a flat connection since it has a non-zero Euler characteristic. It
follows that we can assume $G$ is simple and $(M,\omega)$ is a coadjoint
orbit $(\O, \omega^{\O})$ with its Kirillov--Kostant--Souriau symplectic
structure and with an invariant symplectic connection $\nabla$ with
curvature of Ricci type.

Further, the Euler characteristic of such an orbit is non-zero.
If the vector field $\overline{u}$ were non-zero, then
invariance would imply that it is everywhere non-zero and this
cannot happen. Thus $\overline{u}=0$ and
hence $\nabla$ is locally symmetric ($\nabla R=0$).

Pick a point $\xi_0
\in \O$ and construct a symmetric symplectic triple $(\l,\sigma,\Omega)$
as follows: Let $\a = \{R_{\xi_0}(X,Y) \in \End(T_{\xi_0}\O) \suchthat
X,Y \in T_{\xi_0}\O\}$ and $\l = T_{\xi_0}\O \oplus \a$. The bracket is
defined by
\begin{eqnarray}
[X,Y] &=& R_{\xi_0}(X,Y), \qquad X,Y \in T_{\xi_0}\O;\\{}
[B,X] &=& BX, \qquad B \in \a, X \in T_{\xi_0}\O;\\{}
[B,C] &=& BC-CB, \qquad B,C \in \a,\label{three:eq999}
\end{eqnarray}
$\sigma$ by
\[
\sigma = -\Id_{T_{\xi_0}\O} \oplus \Id_{\a},
\]
and $\Omega$ by
\[
\Omega(X+Y,X'+Y') = \omega_{\xi_0}(X,X'),\qquad X,X' \in T_{\xi_0}, Y,Y'
\in \a.
\]

\begin{lem}
$(\l,\sigma,\Omega)$ is, indeed, a symmetric symplectic triple.
\end{lem}

\begin{pf} There are two things to check to see that $\l$ is a Lie
algebra. Firstly that the brackets defined above belong to $\l$. The
only ones in doubt are the brackets of two elements of $\a$. But $\a$ is
in fact the linear infinitesimal holonomy. This follows since the latter
is spanned by the values of the curvature endomorphism and its covariant
derivatives. The latter vanish by the local symmetry condition.

The second thing to check is the Jacobi identity. Obviously this holds
if all three elements are in $\a$ since this is a Lie algebra. If all three
are in $ T_{\xi_0}\O$ then $[X,[Y,Z]] = - R_{\xi_0}(Y,Z)X$ and the Jacobi
identity is satisfied for these elements by the first Bianchi identity.
When one element is in $ T_{\xi_0}\O$ and two in $\a$ we have
\[
[X,[B,C]]+[B,[C,X]]+[C,[X,B]] = - [B,C]X + BCX - CBX = 0.
\]
Finally, if two elements are in $ T_{\xi_0}\O$ and one in $\a$ we have
\begin{eqnarray*}
[X,[Y,B]] + [Y,[B,X]] + [B,[X,Y]]
&=& -R_{\xi_0}(X,BY) - R_{\xi_0}(BX,Y)\\
&&\qquad\mbox{}+ BR_{\xi_0}(X,Y) - R_{\xi_0}(X,Y)B \\
&=& (B\cdot R_{\xi_0})(X,Y)
\end{eqnarray*}
where $B\cdot R_{\xi_0}$ denotes the natural action of the holonomy Lie
algebra $\a$ on curvature tensors. But $\nabla R=0$ if and only if
$B\cdot R_{\xi_0} = 0, \forall B \in \a$.

The other two properties follow immediately from the definitions.
\end{pf}

If $L$ is the simply connected Lie group associated to $\l$ and $K$
the Lie subgroup associated to the subalgebra $\a$ then $K$ is the
connected component of the fixed point set of the automorphism of $L$
induced by $\sigma$ and $M_1 = L/K$ is a simply connected symmetric space.
$\Omega$ induces a symplectic form $\omega_1$ on $M_1$ which is parallel
for the canonical connection $\nabla_1$.

Consider the point $\xi_0 \in M=\O$ and the point $\xi_1=eK \in M_1$.
There is a linear isomorphism $\phi$ from the tangent space
$T_{\xi_0}M$ to the tangent space $T_{\xi_1}M_1$ so that
\[
\phi\left( R_0(X,Y)Z\right)=R_1(\phi X,\phi Y)\phi Z.
\]
This implies \cite[p.~259, thm.~7.2]{bib:KobNomI} that there
exists an  affine symplectic diffeomorphism
$\psi$ of a neighbourhood $U_0$ of $\xi_0$ in $\O$ onto a  neighbourhood
$U_1$ of $\xi_1$ in $M_1$ such that $\psi_{*_{\xi_0}}=\phi$.

Both $(\O, \omega^{\O}, \nabla)$ and $(M_1,\omega_1,\nabla_1)$ are real
analytic, as is $\psi$, and $\O$ is simply connected whilst $\nabla_1$
is complete. Hence \cite[p.~252, thm.~6.1]{bib:KobNomI} there exists a
unique affine map $\widetilde{\psi} : \O \rightarrow M_1$ such that
$\widetilde{\psi}\vert_{U_0}=\psi$. This map $\widetilde{\psi}$ is
symplectic since it is an analytic extension of  the map $\psi$ which is
symplectic. Symplectic maps are immersions, and hence local
diffeomorphisms when the dimensions are equal as they are in this case.
Hence $\widetilde{\psi}(\O)$ is open. On the other hand, $\O$ is
compact, so $\widetilde{\psi}(\O)$ is compact and thus closed. Hence
$\widetilde{\psi}$ is surjective. It follows that $M_1$ is compact.

{}From the preceding arguments we see that $(M_1,\omega_1,\nabla_1)$ is a
compact simply connected symmetric symplectic space whose curvature is
of Ricci type. The only such space is $\P_n(\C)$ with a multiple of its
standard K\"ahler form $\omega_0$ and the Levi-Civita connection
$\nabla_0$ of the Fubini--Study metric. Since $\O$ and $M_1$ are both
simply connected they are diffeomorphic and hence we have proved Theorem
\ref{thm:1} in the simply connected case.

Next we consider the case where $M$ has a finite fundamental group,
$(M,\omega)$ is $G$-ho\-mo\-ge\-ne\-ous symplectic with a  $G$-invariant
symplectic connection $\nabla$ with curvature of Ricci type. Then the
simply connected covering space $\widetilde{M}$ is compact and carries
such data $\widetilde{\omega},\widetilde{\nabla}$ for the simply
connected covering group $\widetilde{G}$.

It follows that $(\widetilde{M},\widetilde{\omega},\widetilde{\nabla})$
is diffeomorphic to $(\P_n(\C), \omega_0,\nabla_0)$ and hence that $M$
is diffeomorphic to $\P_n(\C)/\Gamma$ where $\Gamma$ is a discrete
subgroup of $PU(n+1)$ acting properly discontinuously on $\P_n(\C)$. But
non-trivial elements of $PU(n+1)$ always have fixed points, so $\Gamma$
must be trivial. This proves Theorem \ref{thm:1}.

\section{} \label{section:4}
We now proceed to give the proof of Theorem \ref{thm:2} indicating along
the way why we restrict ourselves to dimension 4 and why we only obtain a local
result.

Recall that when $(M,\omega)$ is homogeneous and admits a
non-locally-symmetric invariant symplectic connection with
Ricci-type curvature we have the non-zero vector field $\overline{u}$
and  the Ricci endomorphism satisfies
\[
A\overline{u} = 0,
\]
\begin{equation}\label{four:*}
\frac{1}{1+2n}\overline{u} \otimes u - \frac{1+2n}{2(1+n)} A^3 + bA = 0,
\end{equation}
\[
\frac{1+2n}{2(1+n)} A^4 - bA^2 = 0.
\]
Furthermore, if $M$ is compact Lemma \ref{two:lem1} tells us that $b=0$
so that $A^4=0$, and $A^3$ has rank 1:
\[
A^3 = \frac{2(1+n)}{(1+2n)^2}\overline{u}\otimes u.
\]

The $1$-form $u$ is everywhere non-zero so there is a globally defined
vector field $e_1$ with $u(e_1)$ everywhere $\ne 0$.
The vector fields $e_1, e_2=Ae_1, e_3=A^2e_1, e_4=A^3e_1$ form
at each point $x\in M$ a basis of a $4$-dimensional subspace
$V_x$ of the tangent space $T_xM$.
Furthermore, by equation (\ref{four:*})
\[
e_4=\frac{2(1+n)}{(1+2n)^2}u(e_1)\overline{u}.
\]
If we choose the vector field $e_1$ so that $\omega(e_1, e_4)=\epsilon$
with $\epsilon^2=1$, we get $-\frac{2(1+n)}{(1+2n)^2}(u(e_1))^2=\epsilon$
so that $\epsilon=-1$ and $(u(e_1))^2=\frac{(1+2n)^2}{2(1+n)}$ so that
$\overline{u}=u(e_1)e_4$.
 Remark that we can always assume that
$\omega(e_1,e_2)=0$ (by adding to $e_1$ a multiple of $e_3$).
So the symplectic form restricted to $V_x$ writes in the chosen basis
\[
\left(\begin{array}{cccc}
0& 0& 0& -1\\
0& 0& 1& 0\\
0& -1& 0& 0\\
1& 0& 0& 0
\end{array}\right).
\]
The tangent space at each point $x\in M$ writes
\[
T_xM=V_x\oplus V'_x
\]
where $V'_x$ is the $\omega_x$-orthogonal to $V_x$; it is stable under $A$
and, since $A^3$ has rank $1$, $A^3\vert_{V'_x}=0$ but this is not enough
to describe the behaviour of $A$ on $V'_x$.

{}From now on, we restrict ourselves to the $4$-dimensional case.
We define $1$-forms $\alpha, \beta, \gamma, \delta$ such that
\[
\nabla_X e_1 = \alpha(X)e_1 + \beta(X) e_2  + \gamma(X) e_3 + \delta(X) e_4.
\]
Using  formula (\ref{one:1.5}) for $\nabla A$
(i.e. $\nabla_XA = \frac{-1}{2n+1}(X\otimes u + u(e_1)e_4\otimes i(X)\omega)$)
 we obtain
\begin{eqnarray*}
\nabla_X e_2 &=&  \frac{-u(e_1)}{2n+1}  X^1 e_1 +
    \left(\alpha(X)-\frac{u(e_1)}{2n+1}  X^2\right) e_2 \\
              &&\qquad\mbox{}+
    \left(\beta(X)-\frac{u(e_1)}{2n+1}  X^3\right) e_3 +
    \left(\gamma(X)-\frac{2u(e_1)}{2n+1}  X^4\right) e_4,\\
\nabla_X e_3 &=&  \frac{-u(e_1)}{2n+1}  X^1 e_2 +
\left(\alpha(X)-\frac{u(e_1)}{2n+1}  X^2\right) e_3 + \beta(X)e_4,\\
\nabla_X e_4 &=&  \frac{-u(e_1)}{2n+1}  X^1 e_3 +
\left(\alpha(X)-\frac{2u(e_1)}{2n+1}  X^2\right) e_4.
\end{eqnarray*}
On the other hand, formula (\ref{one:1.6}) gives
\[
\nabla_X e_4 = \frac{-u(e_1)}{2n+1}A^2X
\]
so that
\[
\alpha(X) = \frac{u(e_1)}{2n+1} X^2.
\]
The fact that $\nabla$ is symplectic gives the additional
condition that
\[
\gamma(X) = \frac{u(e_1)}{2n+1} X^4.
\]
The connection is thus determined by the two $1$-forms $\beta$ and $\delta$.
The vanishing of the torsion gives the expression of the brackets of the
vector fields $e_j$.

We can now compute the action of the curvature endomorphism on the vector
 fields $e_j$ in two different ways: using the formulas above or
 using the fact that the curvature is of Ricci type.

This yields two identities
\[
d\beta= \frac{3}{2(n+1)}\omega +\frac{u(e_1)}{2n+1}e^2_*\wedge\beta
    +\frac{1}{2(n+1)}e^1_*\wedge e^4_*,
    \]
 \[
 d\delta=2\gamma\wedge \beta -\frac{2}{2(n+1)}e^3_*\wedge e^4_*+
        2\alpha\wedge \delta
 \]
where the $e^j_*$ are $1$-forms so that $e^j_*(e_k)=\delta^j_k$ at
each point.
Using the formulas for the bracket of vector fields we have
\[
de^3_*=\frac{2u(e_1)}{2n+1}e^1_*\wedge e^4_*
-\frac{u(e_1)}{2n+1}e^2_*\wedge e^3_* + e^2_*\wedge \beta,
\]
and substituting $e^2_*\wedge \beta$ in $d\beta$ yields
\[
d(\beta-\frac{u(e_1)}{2n+1}e^3_*)= \frac{2}{n+1}\omega
\]
which is impossible on a compact manifold. This contradiction tells
us that $u$ must vanish and hence that $\nabla$ is locally symmetric.

\medskip

\noindent{\bf Acknowledgement}\ \
The last author is grateful to the Mathematics Department of the
University of Metz for its hospitality during part of this work.

{
}

\end{document}